\documentclass[12pt]{amsart}

\usepackage{amsmath,amssymb,amsthm}
\usepackage{enumerate}
\usepackage[dvips]{color}
\usepackage{version}

\newtheorem{Thm}{Theorem}[section]

\newtheorem{Cor}[Thm]{Corollary}

\newtheorem{Prop}[Thm]{Proposition}

\newtheorem{Conj}[Thm]{Conjecture}
\newtheorem{``Conj"}[Thm]{``Conjecture"}

\theoremstyle{remark}
\newtheorem{Rem}[Thm]{Remark}

\theoremstyle{definition}

\newcommand{\re}{\mathop{\mathrm{Re}}\nolimits}
\newcommand{\im}{\mathop{\mathrm{Im}}\nolimits}

\newcommand{\R}{\ensuremath{\mathbb{R}}}

\newcommand{\Z}{\ensuremath{\mathbb{Z}}}
\newcommand{\Q}{\mathbb{Q}}

\begin{document}

\title[Tropical geometric compactification]
{Collapsing K3 Surfaces and \\ 
Moduli Compactification}
\author{Yuji Odaka}

\author{Yoshiki Oshima}

\maketitle

\begin{abstract}
This note is a summary of our work \cite{OO}, 
which provides an explicit and global moduli-theoretic framework for the collapsing of Ricci-flat K\"ahler metrics and 
we use it to study especially the K3 surfaces case. 
For instance, it allows us to discuss their Gromov-Hausdorff limits along any sequences, which are even not necessarily ``maximally degenerating". 
Our results also give a proof of Kontsevich-Soibelman \cite[Conjecture 1]
{KS} (cf., \cite[Conjecture 6.2]{GW}) in the case of K3 surfaces as a byproduct. 
\end{abstract}


\section{Introduction}

Our paper \cite{OO} 
is a sequel to a 
series  by the first author \cite{TGC.I, TGC.II},  
which compactified both the 
moduli space of compact Riemann surfaces $M_{g} (g\ge 2)$ and that of 
principally polarized abelian varieties $A_{g}$. 
In each case, 
as we actually expect an analogue for any moduli of general polarized 
K\"ahler-Einstein varieties with non-positive scalar curvatures, 
we introduce and study two similar (non-variety) compactifications
 of the moduli space $\mathcal{M}$, which we denote by 
 $\overline{\mathcal{M}}^{\rm GH}$ and 
 $\overline{\mathcal{M}}^{\rm T}$. 
The former $\overline{\mathcal{M}}^{\rm GH}$ is the Gromov-Hausdorff 
compactification with respect to 
\textit{rescaled} K\"ahler-Einstein metrics 
of \textit{fixed diameters} and the latter ``tropical geometric compactification'' 
$\overline{\mathcal{M}}^{\rm T}$ 
should dominate 
the former $\overline{\mathcal{M}}^{\rm GH}$ 
as its boundary $\partial\overline{\mathcal{M}}^{\rm T}$ 
encodes more structure of the 
Gromov-Hausdorff limits (collapses) rather than just distance structure. 
For a precise definition of $\overline{\mathcal{M}}^{\rm GH}$ we employ the same definition as \cite[\S 2.3]{TGC.I}, \cite[\S 2.2]{TGC.II}.\footnote{However, its compactness is unknown at least to the authors in higher dimensional negative scalar curvature case. }
For  $\overline{\mathcal{M}}^{\rm T}$, we have a 
 case by case definition for only particular classes of varieties. 
Here, we recall the structure theorem of 
$\overline{A_{g}}^{\rm GH}$ from \cite[Theorems 2.1 2.3 and Corollary 2.5]{TGC.II}. 

\begin{Thm}[{\cite{TGC.II}}]\label{TGC.Ag.review}
$A_{g}$ can be explicitly compactified as $\overline{A_{g}}^{\rm GH}$ 
whose boundary parametrizes all flat (real) tori $\R^{i}/\Z^{i}$ of diameter  
$1$ where $1\le i\le g$. 
Once we attach the rescaled flat K\"ahler metric
 in the principal polarization with diameter $1$ to each
 abelian variety, the parametrization of metric spaces on
 whole $\overline{A_{g}}^{\rm GH}$ is continuous with
 respect to the Gromov-Hausdorff distance.

\if 
For any algebraic morphism from a punctured curve $f\colon C^{o}=C\setminus \{p\}\to A_{g}$, 
with trivial Raynaud extension, $\lim_{q\to p} f(q)\in \overline{A_{g}}^{\rm T}$ exists and 
the set of possible such limits form a dense subset of $\partial \overline{A_{g}}^{\rm T}$ consists of points with rational coordinates. 
(The Raynaud extension triviality assumption 
 is removed as one consequence of \cite{OO}.) 
 \fi
\end{Thm}

In the above case, 
we simply set $\overline{A_{g}}^{\rm T}:=\overline{A_{g}}^{\rm GH}$. 
On the other hand, in the analogue for $M_{g}$ \cite{TGC.I}, 
we distinguish $\overline{M_{g}}^{\rm GH}$ and $\overline{M_{g}}^{\rm T}$, where 
the boundaries of $\overline{M_{g}}^{\rm GH}$ 
(resp., $\overline{M_{g}}^{\rm T}$)
 parametrize metrized graphs (resp., metrized graphs with integer weights on 
 the vertices). We refer the details to \cite{TGC.I}. 

Our \cite{OO} 
contains the followings. 
\begin{enumerate}
\item \label{MS.Sat.part} We first apply the Morgan-Shalen type compactification  
for 
general Hermitian locally symmetric spaces and identify it with one of the 
Satake compactifications (\cite{Sat1}, \cite{Sat2}). 

\item 
We partially prove that the boundary of the Satake compactification of 
the type which appears in \eqref{MS.Sat.part} parametrizes 
collapses of abelian varieties and 
Ricci-flat K3 surfaces. 
This gives a generalisation of some results in
 \cite{GW}, \cite{Tos}, \cite{GTZ1}, \cite{GTZ2}, \cite{TZ}
 for the K3 surface case. 
For instance, a proof of the conjecture of Kontsevich-Soibelman
\cite[Conjecture 1]{KS} (see also Gross-Wilson \cite[Conjecture 6.2]{GW}),
 which is related to the Strominger-Yau-Zaslow mirror symmetry \cite{SYZ}, for the case of 
K3 surfaces directly follows from our description of collapsing. 
We also give a conjecture for higher dimensional hyperK\"ahler varieties. 
\end{enumerate}
\noindent
Now we move on to a more detailed description. 

\section{General Hermitian symmetric domain}

Let $\mathbb{G}$ be a reductive algebraic group over $\mathbb{Q}$,
 $G=\mathbb{G}(\mathbb{R})$,
 $K$ (one of) its maximal compact subgroup,
 and $D:=G/K$, which we suppose to have a Hermitian symmetric domain structure. 
We moreover assume $D$ is irreducible so that $G$ is simple as a Lie group. 
Suppose that $\Gamma$ is an arithmetic subgroup of $\mathbb{G}(\mathbb{Q})$,
 which acts on $D$. 
Hence we can discuss Hermitian locally symmetric space $\Gamma\backslash D$. 

Satake \cite{Sat1}, \cite{Sat2} 
constructed compactifications of Riemannian locally symmetric spaces $G/K$ 
associated to irreducible projective representations $\tau\colon G\to
PGL(\mathbb{C})$ satisfying certain conditions. They are 
 stratified as: 
\[\overline{\Gamma\backslash D}^{\rm Sat, \tau}
 =\Gamma\backslash D\sqcup \bigsqcup_{P}
(\Gamma\cap Q(P))\backslash M_{P}/(K\cap M_{P}).\] 
Here, $P$ runs over all the $\mu(\tau)$-connected rational parabolic subgroups, 
 $P=N_{P}A_{P}M_{P}$ denotes the Langlands decomposition, and 
 $Q(P)$ is the $\mu(\tau)$-saturation of $P$. 
We are particularly interested in the case when $\tau$ is the adjoint
 representation $\tau_{\rm ad}$. 
 
 On the other hand, 
 given any toroidal compactification \cite{AMRT} 
 for $\Gamma\backslash D$, we can apply the Morgan-Shalen type compactification 
 to it as \cite[Appendix]{TGC.II} (following \cite{MS, BJ16}). 
The Morgan-Shalen type compactification
 $\overline{\Gamma\backslash D}^{\rm MSBJ}$ obtained in this way
 is independent of the cone decomposition 
 for the toroidal compactification  \cite[A.13, A.14]{TGC.II}. 

We now compare these two compactifications.

\begin{Thm}\label{MS.Sat}
Let  $\Gamma\backslash D$ be a locally Hermitian symmetric space. 
Consider its toroidal compactification 
 and the associated (generalised) Morgan-Shalen compactification
 $\overline{\Gamma\backslash D}^{\rm MSBJ}$. 
Then this is homeomorphic to the Satake compactification 
$\overline{(\Gamma\backslash D)}^{\rm Sat,\tau_{\rm ad}}$ 
for the adjoint representation $\tau_{\rm ad}$ of $G$.
\end{Thm}

In the following we make an ``elementary'' but important observation
 on a rationality phenomenon of the limits along one parameter holomorphic 
family, which we expect to fit well with the recent approach to extend 
the theta functions in \cite{GS.JDG} etc.

\begin{Prop}\label{MS.lim}
Suppose $U\subset \overline{U}^{\rm hyb}(\mathcal{X})$ is a Morgan-Shalen-Boucksom-Jonsson compactification associated to an arbitrary dlt stacky pair $(\mathcal{X},\mathcal{D})$ of boundary coefficients $1$ 
(\cite{TGC.II}) with $\mathcal{U}:=\mathcal{X}\setminus \mathcal{D}$, 
its coarse moduli space $\mathcal{U}\to U$. Then for any holomorphic morphism 
$\Delta^{*}:=\{z\in \mathbb{C}\mid 0<|z|<1\}\to \mathcal{U}$ which extend to 
$\Delta:=\{z\in \mathbb{C}\mid |z|<1\}\to \mathcal{X}$, it induces a continuous map  
$\Delta\to \overline{U}^{\rm hyb}(\mathcal{X})$, i.e., the limit exists. 
Furthermore, such possible limits in $\Delta(\mathcal{D})$ are 
characterized as points with rational coordinates. 
\end{Prop}

\begin{Cor}[corollary to Theorem~\ref{MS.Sat} and Proposition~\ref{MS.lim}]\label{}
Take an arbitrary holomorphic map $f\colon \Delta^*\to \Gamma\backslash D$, 
which extends to a map to a toroidal compactification of $\Gamma\backslash D$. 
Then $f$ also extends to a map $\Delta\to 
\overline{\Gamma\backslash D}^{\rm Sat,\tau_{ad}}$ 
where $0$ is sent to a point with rational coordinates, i.e., 
a point in the dense subset
 $(C(F)\cap U(F)\otimes \Q)/\Q_{>0}\subset C(F)/\mathbb{R}_{>0}$. 
\end{Cor}
This is partially proved 
in the case of $A_{g}$ 
in \cite{TGC.II} by using degeneration data in 
\cite{FC90}. 
\begin{Rem}
Although we assume that $G$ is simple in this section,
 our Morgan-Shalen type compactification
 construction \cite[Appendix]{TGC.II} still works for non-simple $G$. 
Thus, our construction also gives a new Satake-type compactification
 for non-simple $G$, e.g., of the Hilbert modular varieties. 
\end{Rem}


\section{Abelian varieties case} 

We identify our tropical geometric compactification $\overline{A_g}^{\rm T}$ (\cite{TGC.II})
of $A_g$ with the adjoint type Satake compactification. 

\begin{Thm}\label{Ag.TGC.Satake.MS}
There are canonical homeomorphisms between the three compactifications  
\[\overline{A_g}^{\rm T} \cong \overline{A_g}^{\rm Sat,\tau_{\rm ad}}\cong  \overline{A_g}^{\rm MSBJ},\]
extending the identity on $A_g$. 
\end{Thm}
The second canonical homeomorphism is a special case of Theorem \ref{MS.Sat} and 
the first is essentially reduced to matrix computations. 


In \cite{OO}, we also give a purely moduli-theoritic 
reexplanation 
of the structure theory of one parameter degenerations of 
abelian 
varieties in \cite{Mum72.AV}, \cite{FC90}, 
after the above Theorem \ref{Ag.TGC.Satake.MS} as follows. 

\begin{Thm}
Take a holomorphic maximally degenerating family of principally polarized 
abelian varieties 
$\pi\colon (\mathcal{X},\mathcal{L})\to \Delta$. 
Consider the rescaled Gromov-Hausdorff limit $B(\mathcal{X},\mathcal{L})$ of 
diameter $1$ as in Theorem \ref{TGC.Ag.review} (\cite{TGC.II}) 
and its discrete Legendre transform 
$\check{B}(\mathcal{X},\mathcal{L})$ (\cite{GS11}, \cite{KS}). 

Then we can enhance the underlying integral affine structure of 
$\check{B}(\mathcal{X},\mathcal{L})$ 
as $K$-affine structure (in the sense of \cite[\S 7.1]{KS}) naturally 
via the data of $\pi$. 
Furthermore, such $K$-affine structure recovers $\pi$ 
up to an equivalence relation generated by 
base change (replace $t$ by $t^{a}$ with $a\in \Q_{>0}$). 
\end{Thm}



\section{Moduli of Algebraic K3 surfaces}
\subsection{Satake compactification}
\label{K3.Sat.sec}

Let $\mathcal{F}_{2d}$ be the moduli space of 
 polarized K3 surfaces of degree $2d$ possibly with ADE singularities. 
Its structure is known as follows. Let 
 $\Lambda_{\rm K3}:=E_{8}(-1)^{\oplus 2}\oplus U^{\oplus 3}$ be the K3 lattice
 and fix a primitive vector $\lambda_{2d}$ with $(\lambda_{2d},\lambda_{2d})=2d$
 and $\Lambda_{2d}:=\lambda_{2d}^{\perp}$.
The complex manifold
$$\Omega(\Lambda_{2d}):=\{[w]\in \mathbb{P}(\Lambda_{2d}\otimes \mathbb{C})\mid
 (w,w)=0,\ (w,\bar{w})>0\}.$$ 
 has two connected components.
We choose one component and denote by $\mathcal{D}_{\Lambda_{2d}}$. 
Let $O(\Lambda_{\rm K3})$ denote the isomorphism group
 of the lattice $\Lambda_{\rm K3}$
 preserving the bilinear form and  
set
\begin{align*}
\tilde{O}(\Lambda_{2d}):=\{g|_{\Lambda_{2d}} : g\in O(\Lambda_{\rm K3}),\, 
g(\lambda_{2d})=\lambda_{2d}\}.
\end{align*}
The group $\tilde{O}(\Lambda_{2d})$ naturally acts on $\Omega(\Lambda_{2d})$.
We define $\tilde{O}^{+}(\Lambda_{2d})$ to be 
 the index two subgroup of $\tilde{O}(\Lambda_{2d})$
 consisting of the elements preserving each connected component
 of $\Omega(\Lambda_{2d})$. 
Then it is well-known that 
\begin{align*}
\mathcal{F}_{2d}
\simeq 
\tilde{O}^{+}(\Lambda_{2d})\backslash \mathcal{D}_{\Lambda_{2d}}
\simeq
\tilde{O}(\Lambda_{2d})\backslash \Omega(\Lambda_{2d}).
\end{align*}
Let $\overline{\mathcal{F}_{2d}}^{{\rm Sat},\tau_{\rm ad}}$
 (or simply $\overline{\mathcal{F}_{2d}}^{{\rm Sat}}$ in our papers)
 be the Satake compactification of  $\mathcal{F}_{2d}$
 corresponding to the adjoint representation of $O(2,19)$.
It decomposes as
\[\overline{\mathcal{F}_{2d}}^{{\rm Sat}}=
\mathcal{F}_{2d}\sqcup \bigcup_{l} \mathcal{F}_{2d}(l)
\sqcup \bigcup_{p} \mathcal{F}_{2d}(p),\]
where $l$ runs over one-dimensional isotropic
 subspaces of $\Lambda_{2d}\otimes \mathbb{Q}$,
 and $p$ runs over two-dimensional isotropic subspaces of
 $\Lambda_{2d}\otimes \mathbb{Q}$. 
Also, we simply define the tropical geometric compactification of 
$\mathcal{F}_{2d}$ as this $\overline{\mathcal{F}_{2d}}^{\rm Sat}$. 
The boundary component $\mathcal{F}_{2d}(l)$ is given as 
\[
\mathcal{F}_{2d}(l)
= \{v\in (l^{\perp}/l) \otimes \mathbb{R} \mid (v,v)>0\}/\sim.
\]
Here $v \sim v'$ if $g\cdot v=c v'$ for some $g\in \tilde{O}^{+}(\Lambda_{2d})$
 and $c\in \mathbb{R}^{\times}$.
We have $\mathcal{F}_{2d}(l)=\mathcal{F}_{2d}(l')$
 if $g\cdot l=l'$ for some $g\in \tilde{O}^{+}(\Lambda_{2d})$
 and $\mathcal{F}_{2d}(l)\cap\mathcal{F}_{2d}(l')=\emptyset$ if otherwise.
Since $(l^{\perp}/l) \otimes \mathbb{R}$ has signature $(1,18)$,
 there is an isomorphism
\begin{multline*}
\{v\in (l^{\perp}/l) \otimes \mathbb{R} \mid (v,v)>0\}/\mathbb{R}^{\times}\\
\simeq O(1,18)/ O(1)\times O(18)
\end{multline*}
and hence $\mathcal{F}_{2d}(l)$
 is an arithmetic quotient of $O(1,18)/O(1)\times O(18)$.
The other component $\mathcal{F}_{2d}(p)$ is a point
 and $\mathcal{F}_{2d}(p)=\mathcal{F}_{2d}(p')$
 if and only if $g\cdot p=p'$ 
 for some $g\in \tilde{O}^+(\Lambda_{2d})$. 
Therefore, if we take representatives of $l$ and $p$
 from each equivalence class, we get a finite decomposition: 
\[\overline{\mathcal{F}_{2d}}^{{\rm Sat}}=
\mathcal{F}_{2d}\sqcup \bigsqcup_{l} \mathcal{F}_{2d}(l)
\sqcup \bigsqcup_{p} \mathcal{F}_{2d}(p).\]

\subsection{Tropical K3 surfaces}\label{trop.K3.1}

In our paper, what we mean by \textit{tropical polarized K3 surface} 
is a topological space $B$ 
homeomorphic to the sphere $S^{2}$, with an 
affine structure away from certain finite points ${\rm Sing}(B)$, 
with a metric which is Mong\'e-Ampere metric $g$ 
with respect to the affine structure on $B\setminus {\rm Sing}(B)$. 
Studies of such object as tropical version of 
K3 surfaces are pioneered in well-known papers of Gross-Wilson \cite{GW} and 
Kontsevich-Soibelman \cite{KS}. 

Here we assign such tropical K3 surface to each point in
 the boundary component $\mathcal{F}_{2d}(l)$ as follows. 
Let $l$ be an oriented one-dimensional isotropic
 subspace of $\Lambda_{2d}\otimes \mathbb{Q}$.
Write $e$ for the primitive element of $l$ such that
 $\mathbb{R}_{> 0} e$ agrees with the orientation of $l$.
Take a vector $v\in (l^{\perp}/l)\otimes \mathbb{R}$ such that $(v,v)>0$.
Write $[e,v]$ for the corresponding point in $\mathcal{F}_{2d}(l)$.
Then there exists a (not necessarily projective) K3 surface $X$
 and a marking 
 $\alpha_{X}\colon H^2(X,\mathbb{Z})\to \Lambda$ with 
\begin{itemize}
\item $\alpha_X(H^{2,0})\subset \mathbb{R}\lambda+\sqrt{-1}\mathbb{R}v$,
\item $\alpha_{X}^{-1}(e)$ is in the closure of 
K\"ahler cone.
\end{itemize}
The pair $(X,\alpha_{X})$ is unique up to isomorphisms.

Let $L$ be a line bundle on $X$ such that $\alpha_{X}([L])=e$.
Then we get an elliptic fibration $f:X\to B(\simeq \mathbb{P}^1)$.
Take a holomorphic volume form $\Omega$ on $X$ such that
 $\alpha_X([\re \Omega])=\lambda$.
The map $f$ is a Lagrangian fibration with respect to
 the symplectic form $\re \Omega$.
Hence it gives an affine manifold structure
 on $B\setminus\Delta$,
 where $\Delta$ denotes the finite set of singular points.
Similarly, the imaginary part $\im \Omega$ gives another affine manifold
 structure on $B\setminus\Delta$.

We endow the base space $B$ with 
the McLean metric on the base $B$ (\cite{ML}), 
where we regard $f$ as special Lagrangian fibration after hyperK\"ahler rotation. 
A straightforward calculation shows that this coincides with the 
``special K\"ahler metric" $g_{\it sp}$ introduced and studied 
in \cite{DW96, Hit,Freed99} and appears as the metric on $\mathbb{P}^1$ 
in \cite{GTZ2}. We rescale the metric to make its diameter $1$
 and denote this obtained tropical K3 surface by $\Phi_{\rm alg}([e,v])$. 
\begin{Rem}
Recall the concepts of the \textit{class of metric} (\textit{metric class}) 
and the \textit{radiance obstruction} 
 of Mong\'e-Amp\'ere manifolds $B$ with singularities. They are 
introduced in \cite{KS} and discussed in \cite{GS.logI} in more details. 
We denote them by 
$k(B)\in H^1(B, i_*\tilde{\Lambda}^{\vee}\otimes \mathbb{R})$ and 
$c(B)\in H^1(B,i_*\Lambda)$, respectively.
Here, $\Lambda$ is the 
affine structure as a $\mathbb{Z}^{\it dim(B)}$-local system in tangent 
bundle $T(B\setminus \Delta)$, $-^{\vee}$ denotes $-$'s 
dual local system, $\tilde{\Lambda}^{\vee}$ is local system of affine 
functions. In particular, we naturally have a morphism of local systems 
$f\colon \tilde{\Lambda}^{\vee}\to \Lambda^{\vee}$ which induces 
$f_*\colon H^1(B, i_*\tilde{\Lambda}^{\vee})\to H^1(B, i_*\Lambda^{\vee})$. 
It is also easy to see that, if we slightly change the definition of the metric class, 
to extract its ``linear" part as $f_* k(B)$. 
Then, 
it naturally recovers the data $\overline{v}\in (e^{\perp}\otimes \mathbb{R}/\mathbb{R}e)$ 
i.e., we have 
$f_*k(\Phi_{\rm alg}([e,v]))=[v],$ 
under the natural identification 
$H^1(\Phi_{\rm alg}([e,v]), i_*\Lambda^{\vee}\otimes \mathbb{R})\hookrightarrow
(e^{\perp}\otimes \mathbb{R}/\mathbb{R} e)$ 
which comes from the Leray spectral sequence applied to the elliptic fibration 
$X\twoheadrightarrow \Phi_{\rm alg}([e,v])$ in \S \ref{trop.K3.1}. 
Our results in \cite{TGC.II} and Theorem~\ref{Ag.TGC.Satake.MS} for $A_{g}$ 
can be re-interpretted similarly (but with weight $1$). 
\end{Rem}

\begin{Rem}
\textit{Yuto Yamamoto} \cite{Yam} has some ongoing interesting work 
which seems to be related to our works, 
where he constructs a sphere with an integral affine structure 
from the tropicalization of an anticanonical hypersurface 
in a toric Fano 3-fold, and computes its radiance obstruction. 
\end{Rem}

\subsection{Gromov-Hausdorff collapse of K3 surfaces}
For a point in $\mathcal{F}_{2d}$
 we have a corresponding polarized K3 surface $(X,L)$,
 equipped with a natural Ricci-flat metric.
For $[e,v]\in \mathcal{F}_{2d}(l)$ we defined in a previous section
$\Phi_{\rm alg}([e,v])$.  
For a point in $\mathcal{F}_{2d}(p)$
 we assign a (one-dimensional) segment,
 which we denote by $\Phi_{\rm alg}(\mathcal{F}_{2d}(p))$. 
Let us normalize these metric spaces so that their diameters are one.
We thus obtained a map
 $\Phi_{\rm alg}\colon \overline{\mathcal{F}_{2d}}^{{\rm Sat}}
 \to \{\text{compact metric spaces with diameter one}\}$. 
 Here, we associate Gromov-Hausdorff distance to the right hand side (target space) 
 and denote it by ${\it CMet}_{1}$.

\begin{Conj}\label{K3.Main.conjecture}
The map 
\[\Phi_{\rm alg}\colon \overline{\mathcal{F}_{2d}}^{{\rm Sat}}
 \to {\it CMet}_{1}\]
given above is continuous.
\end{Conj}

We would like to simply set 
the tropical geometric compactification of $\mathcal{F}_{2d}$ as 
$\overline{\mathcal{F}_{2d}}^{{\rm T}}:=\overline{\mathcal{F}_{2d}}^{{\rm Sat}}$. 
Indeed, if Conjecture~\ref{K3.Main.conjecture} holds, 
we get a continuous map $\overline{\mathcal{F}_{2d}}^{{\rm Sat}}\to 
\overline{\mathcal{F}_{2d}}^{\rm GH}$ and we also observe that 
each $\mathcal{F}_{2d}(l)$ encodes affine structure of the limit tropical K3 
surface as well. (This answers a question of Prof.\ B.~Siebert in 2016 to the first author, 
regarding if one can associate tropical affine structure to limit of any 
collapsing \textit{sequence}). 
So far, we have partially confirmed the conjecture. 
The case of ($A_{1}$-singular flat) Kummer surfaces, with $3$-dimensional moduli,
 are easily reduced to \cite{TGC.II}.
More generally, we have proved the following. 
In particular, the conjecture \ref{K3.Main.conjecture} holds 
at least away from finite points. 
\begin{Thm}\label{K3.Main.Conjecture.18.ok}The map 
$\Phi_{\rm alg}$ is continuous on $\overline{\mathcal{F}_{2d}}^{\rm Sat}\setminus (\bigcup_{p}\mathcal{F}_{2d}(p))$. 
It is continuous also 
when restricted to the boundary $\partial{\overline{\mathcal{F}_{2d}}^{\rm Sat}}
=\overline{\mathcal{F}_{2d}}^{\rm Sat}\setminus \mathcal{F}_{2d}$. 
\end{Thm}
The proof of the former half of the statements 
involves some symmetric space theory,
hyperK\"ahler geometry, 
algebraic geometry of moduli, and a priori analytic estimates. 
The estimates heavily 
 depends on \cite{Tos,GW,GTZ1,GTZ2,TZ} 
 and their extensions. 
One nontrivial part of the extension
 is, for instance, to make many of the $C^{2}$-estimations in \textit{op.cit}
 following methods of \cite{Yau} 
 locally uniform with respect to a family of elliptic K3 surfaces
 even along degenerations to orbifolds. 

During our work,  we learnt that 
\textit{Kenji Hashimoto, Yuichi Nohara, Kazushi Ueda} \cite{HNU} 
also studied the Gromov-Hausdorff collapses along certain $2$-dimensional 
subvariety of $\mathcal{F}_{2d}$, i.e., the moduli 
of $E_{8}^{\oplus 2}\oplus U(\oplus \langle -2 \rangle)$-polarized K3 surfaces.
Moereover, a result of Hashimoto and Ueda \cite{HU} implies that
 the restriction of $\Phi_{\rm alg}$ to the boundary
 is a generically two-to-one map. 
We appreciate their gentle discussion with us. 

Theorem \ref{K3.Main.Conjecture.18.ok} (resp., Conjecture \ref{K3.Main.conjecture}) combined with 
Proposition~\ref{MS.lim} determines the Gromov-Hausdorff limits of Type III (resp., Type II) one parameter family of 
Ricci-flat algebraic K3 surfaces, which solves a conjecture of 
Kontsevich-Soibelman \cite[Conjecture 1]{KS}, Todorov,  and Gross-Wilson (cf., e.g., \cite[Conjecture 6.2]{Gross}) 
 in the K3 surfaces case. 

In the next section, we discuss collapsing of 
general K\"ahler K3 surfaces, which are not necessarily algebraic. 


\section{Moduli of K\"ahler K3 surfaces}
It is known 
(cf., \cite{Tod}, \cite{Looi}, \cite{KT})
that the moduli space
 of all Einstein metrics on a K\"ahler K3 surfaces 
 (including orbifold-metrics)
 has again a structure of the locally Riemannian symmetric space: 
$$O(\Lambda_{\rm K3})\backslash SO_{0}(3,19)/
(SO(3)\times SO(19)),$$ 
which we denote by $\mathcal{M}_{\rm K3}$.  
An enriched version encoding also complex structures of the K3 surfaces is 
$$\R_{>0}\times (O(\Lambda_{\rm K3})\backslash SO_{0}(3,19)/
(SO(2)\times SO(19))).$$
Roughly speaking, 
this is a union of K\"ahler cones of ADE K3 surfaces with marking 
of the minimal resolutions. 

Thus we can again compare  a 
Satake compactification of $\mathcal{M}_{\rm K3}$ 
with the Gromov-Hausdorff compactification. 
Inside the Satake compactification for the adjoint representation, 
we consider an open locus (a partial compactification of $\mathcal{M}_{\rm K3}$) 
$\mathcal{M}_{\rm K3} \sqcup 
\mathcal{M}_{\rm K3}(a), $
where $\mathcal{M}_{\rm K3}(a)$ denotes
 the $36$-dimensional boundary stratum
 corresponding to an isotropic rational line
 $l=\Q e$ in $\Lambda_{\rm K3}\otimes \Q$,
 with primitive integral generator $e$,
 which are unique up to $O(\Lambda_{\rm K3})$. 
Then for each point $p=\langle e, v_{1},v_2 \rangle$ in strata 
$\mathcal{M}_{\rm K3}(a)$, we consider 
the marked (possibly ADE) K3 surface $X_{p}$ with period $\langle v_{1},v_2\rangle$. Then it is known that there is an elliptic K3 surface structure on 
$X_p$ with the fiber class $e$. Then we define 
 $\Phi(p)$ as its base biholomorphic to  $\mathbb{P}^1$ with the McLean metric,
 which only depends on $\langle v_1,v_2\rangle$. 
Similarly to the projective case Theorem~\ref{K3.Main.Conjecture.18.ok},
 \cite{OO} proves that for non-algebraic situation: 

\begin{Thm}\label{K3.Main.conjecture2}
The map 
$$
\Phi\colon 
\mathcal{M}_{\rm K3} \sqcup 
\mathcal{M}_{\rm K3}(a) \to {\it CMet}_{1}$$ 
given above is continuous. 
Here, we put the Gromov-Hausdorff topology for the right hand side. 
\end{Thm}

In \cite{OO}, we further explicitly 
define an extension to the whole Satake compactification 
 $\Phi\colon \overline{\mathcal{M}_{\rm K3}}^{\rm Sat}\to {\it CMet}_{1}$, 
and conjecture that this is still continuous with respect to the 
Gromov-Hausdorff topology. 
For the boundary strata other than $\mathcal{M}_{\rm K3}(a)$,
 we assign flat tori $\R^i/\Z^i\ (i=1,2,3)$
 modulo $(-1)$-multiplication. We show that $\Phi$ restricted to the closure of the locus which 
 parametrizes $\R^4/\Z^4$ modulo $\pm 1$, that includes those boundary strata, is continuous. 
Furthermore, 
 we also prove the restriction of $\Phi$ to
 the closure of $\mathcal{M}_{\rm K3}(a)$ is continuous
 by using Weierstrass models. 


\section{Higher dimensional case}
We expect that our results for K3 surfaces naturally extend to  
higher dimensional compact hyperK\"ahler manifolds. Let us focus on 
algebraic case in this notes. We set up as follows. 
Fix any connected moduli $M$ of polarized $2n$-dimensional irreducible holomorphic 
symplectic manifolds $(X,L)$ whose second cohomology 
$H^{2}(X,\Z)$ is isomorphic (as a lattice) 
to $\Lambda$. By 
\cite{Ver, Mark} (\cite[Theorem 3.7]{GHS}), 
it is a Zariski open subset of a Hermitian locally symmetric space of orthogonal type $\Gamma\backslash \mathcal{D}_{M}$. 

Then (a rough version of) our 
conjecture for algebraic case (in \cite{OO}) is as follows: 

\begin{Conj}
There 
is a continuous map $\Psi$ (call ``geometric realization map") 
from the Satake compactification $(M\subset)\overline{\Gamma\backslash \mathcal{D}_M}^{\rm Sat,\tau_{ad}}$ 
with respect to the adjoint representation 
to the Gromov-Hausdorff compactification of $M$, extending the identity map on $M$. 
The $(b_{2}(X)-4)$-dimensional boundary strata of $\overline{\Gamma\backslash \mathcal{D}_M}^{\rm Sat,\tau_{ad}}$  parametrize via $\Psi$ the projective space 
$\mathbb{P}^{n}$ with special K\"ahler metrics in the sense of \cite{Freed99}
and the 
metric space parametrized by 
$0$-dimensional cusps are all homeomorphic to the closed ball of dimension $n$. 
\end{Conj}
At the moment of writing this notes, 
the authors have only succeeded in proving that $(M\subset )\Gamma\backslash \mathcal{D}_M$ 
is the moduli of polarized symplectic varieties with continuous (non-collapsing) weak Ricci-flat K\"ahler metrics,  
and making some progress on the necessary algebro-geometric preparations in particular for 
the case of K3$^{[n]}$-type. 

\begin{Rem}[Calabi-Yau case]
In \cite{OO}, we also propose an extension 
of Conjecture~\ref{K3.Main.conjecture} for general Calabi-Yau varieties under some technical conditions, 
although there are much fewer evidences in that case. 
\end{Rem}


\textbf{Acknowledgement}
We appreciate for giving us the chances to talk on \cite{OO} 
in various countries and cities. The first was at a  
talk by the first author at a Clay conference held at Oxford in September 
2016, when 
Theorems~\ref{K3.Main.Conjecture.18.ok} and \ref{K3.Main.conjecture2} 
were only partially proved and claimed, 
whose confirmation in the form of this notes has taken long time. 
In particular, we appreciate 
Kenji Hashimoto, Shouhei Honda, Radu Laza, Daisuke Matsushita, Shigeru Mukai, Yoshinori Namikawa, 
Bernd Siebert, Cristiano Spotti, Song Sun, 
Yuichi Nohara, Kazushi Ueda, Ken-ichi Yoshikawa 
for helpful discussions. There are plans of some lecture series by the first author on this topic during the next fall semester in Nagoya, Tokyo. 
The first author is partially supported by JSPS Grant-in-Aid (S), No. 16H06335, 
Grand-in-Aid for Early-Career Scientists No. 18K13389. 
The second author is partially supported by JSPS KAKENHI Grant No.\ 16K17562.


\vspace{5mm} \footnotesize \noindent
Contact (Yuji Odaka): {\tt yodaka@math.kyoto-u.ac.jp} \\
Department of Mathematics, Kyoto University, Kyoto 606-8285. JAPAN \\

\vspace{5mm} \footnotesize \noindent
Contact (Yoshiki Oshima): {\tt oshima@ist.osaka-u.ac.jp} \\
Graduate School of Information Science and Technology,
 Osaka University, Suita, Osaka 565-0871, JAPAN \\

\end{document}